\documentclass[11pt,english,runninghead]{article}
\usepackage[T1]{fontenc}
\usepackage[latin9]{inputenc}
\usepackage[letterpaper]{geometry}
\usepackage{amsthm}
\usepackage{amsmath}
\usepackage{graphicx}
\usepackage{setspace}
\usepackage{amssymb}
\usepackage{esint}

\makeatletter

\numberwithin{equation}{section}
  \theoremstyle{plain}
   
  \theoremstyle{remark}
  
  \theoremstyle{plain}

\usepackage{mathrsfs}

\makeatother

\usepackage{babel}

\makeatother

\usepackage{babel}

\begin{document}

\Large

\begin{center}
{\bf  Wellposedness of a  swimming model  in the 3-$D$ incompressible fluid governed 
by the  nonstationary Stokes equation.}
\end{center}

\normalsize

\bigskip

\begin{center}

A.Y. Khapalov \footnote{This work was supported in part by NSF
Grant DMS-1007981.} \\ Department of Mathematics\\Washington State
University,  Pullman, WA 99164-3113; \\e-mail: khapala@math.wsu.edu \\

\end{center}
\begin{abstract}
We introduce and investigate the wellposedness of a model describing the self-propelled  motion of a {\em small} abstract {\em swimmer} 
in the  3-$D$ incompressible fluid governed 
by the  nonstationary Stokes equation, typically associated with the low Reynolds numbers.  It is  assumed that the swimmer's body  consists of 
finitely many subsequently connected parts, identified with the fluid they occupy,  linked by
 rotational  and elastic Hooke's forces. In this paper we are attempting to extend the 2-$D$ version of this model, introduced in \cite{Kh11}-\cite{Kh2}, to the 3-$D$ case.
Models like this are of interest in biological
and engineering applications dealing with the study and design of propulsion systems in
fluids. 
\end{abstract}

\bigskip
 {\bf Key words:}
Swimming models, hybrid systems,  nonstationary Stokes fluid.

\bigskip
{\bf 1. Introduction and problem formulation.}  It seems that the first quantitative research  in the area of swimming phenomenon was
aimed at the biomechanics of specific  biological species: Gray
 \cite{Gra1}(1932), Gray and Hancock  \cite{Gra2} (1951), Taylor  \cite{Tay1} (1951),  \cite{Tay2} (1952),  Wu \cite{Wu} (1971), Lighthill  \cite{Lig} (1975),
and others. These efforts  resulted in the derivation of a  number of mathematical models (linked the size of Reynolds number) 
for swimming motion in 
the whole $R^2$- or $ R^3$-spaces with the  swimmer to be used as 
the reference frame, see, e.g.,
Childress  \cite{Chi}  (1981) and the references therein. Such approach however
requires  some modification if one wants   to  track  the actual position of swimmer in
a fluid.

A different modeling approach was proposed by Peskin  in the computational
mathematical biology  (see Peskin  \cite{Pes1}  (1975), Fauci and Peskin  \cite{Fau1} (1988), Fauci  \cite{Fau2}  (1993), Peskin and McQueen  \cite{Pes2}  (1994) and the references therein), where a 
swimmer is modeled as an
immaterial  {\em immersed boundary} identified with the fluid, further discretized for computational purposes on some grid.
 In this case a fluid
equation is to be complemented by a coupled {\em infinite} dimensional
differential equation for the aforementioned ``immersed boundary''. 

{\em In this paper}  we intend  to deal with the swimming phenomenon in the framework  of {\em non-stationary} PDE's along the {\em immersed body} approach summarized in   Khapalov \cite{Kh2}  (2005-2010).
Namely, in \cite{Kh11} (2005), inspired by  the ideas of the above-cited Peskin's  method,   we introduced a   2-$D$  model for {\em ``small'' flexible swimmers} assuming that their  bodies are identified with the fluid  occupying their shapes. This approach views  such a swimmer as the already {\em discretized} aforementioned immersed boundary  supported on the respective grid cells, see, e.g.,  Fig. 1 below.
Our model offered  two novel features:   (1) it was set in  {\em  a bounded} domain with  (2) governing equations to be  a fluid  equation  coupled with  a {\em system of ODE's} describing the spatial position of swimmer within the  space domain.   We established the wellposedness of this model  up to the contact of  a swimmer at hand either with the boundary of space domain  or with itself. 
The need of such type of  models was motivated by the intention  to investigate
controllability properties of swimming phenomenon (see \cite{Kh2}). 

{\em Our goal} in this paper is to investigate the wellposedness of a 3-$D$ version of this model.

\begin{center}
\includegraphics[scale=0.6]{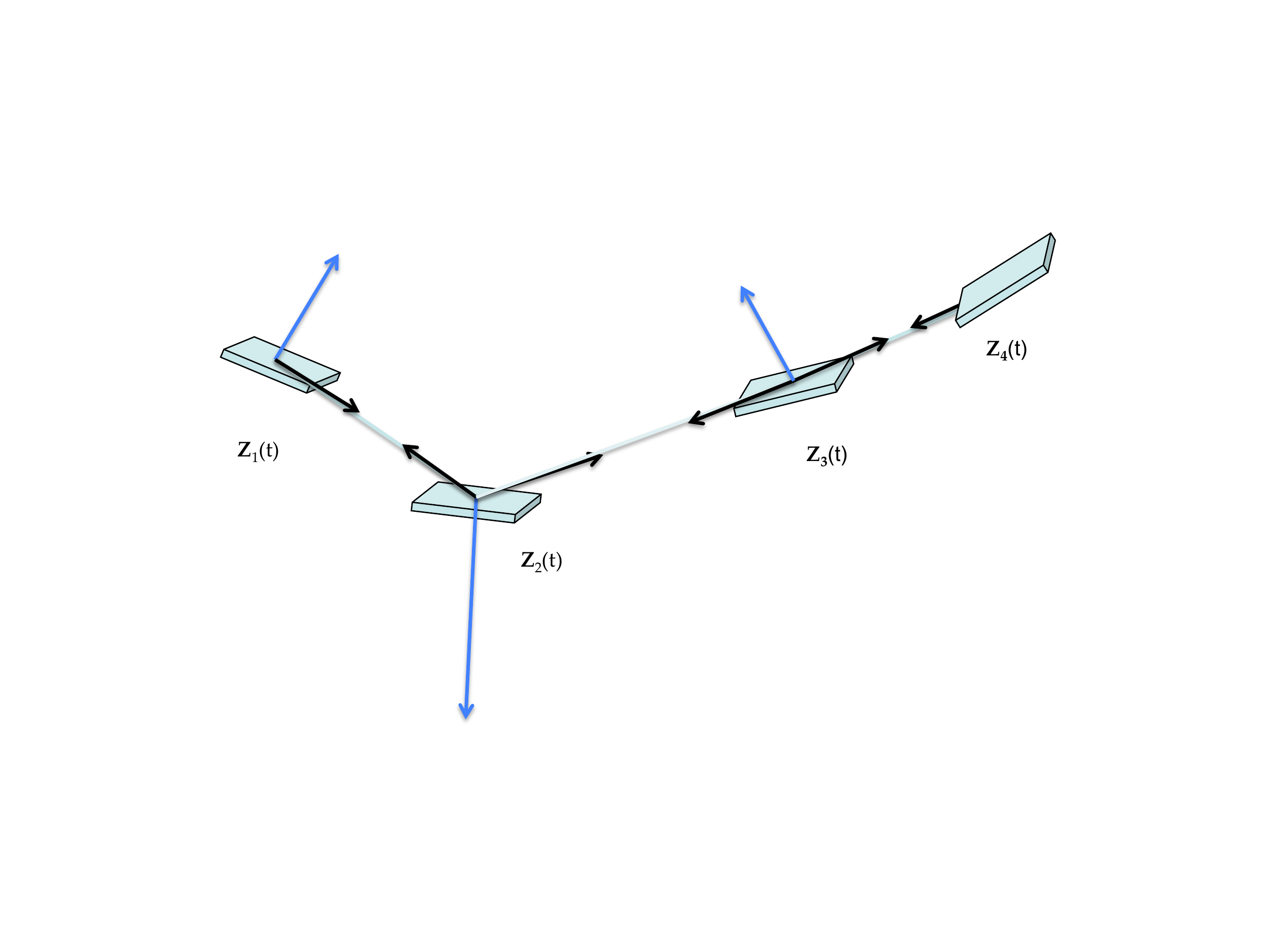} 
\par\end{center}

{\sf Figure 1: 4-parallelepipeds swimmer with all elastic forces and   rotation forces about $ z_2(t)$ only.}

\bigskip

{\bf Further remarks on bibliography.}  It should be noted that,   the
classical mathematical issue of  wellposedness of a swimming model  as a system of PDE's 
  for the first time  was apparently addressed by
Galdi  \cite{Gal}  (1999) for a model of swimming micromotions  in $ R^3$
(with the swimmer serving  as the reference frame).

Another available approach to modeling of swimming motion (apparently,  initiated by the work  Shapere and Wilczeck  \cite{Sha} (1989)) exploits the idea that the swimmer's body shape transformations during the actual swimming process  can be viewed as  a set-valued map in time. The respective models describe swimmer's position via such maps, see  \cite{Cur} (1981),   \cite{San2} (2008), \cite{Dal} (2011)  and the references therein.  Some models treat these maps as a priori prescribed, in which case the crux of the problem is to identify which maps  are admissible, i.e., compatible with the principle of self-propulsion  of swimming locomotion. In the case when the aforementioned motion map is not a priori prescribed (i.e.,  it will be defined at each moment of time by swimmer's internal forces and the interaction of its body with the resisting surrounding medium), the model  will have to include extra equations, see, e.g.,  \cite{Fau3} in the framework of the immersed boundary method and the references therein.

More recently, a number of significant
efforts, both theoretical and experimental, were made to study models of possible bio-mimetic mechanical devises
which employ the change of
their geometry, inflicted by internal forces, as the means for self-propulsion, see, e.g.,
S. Hirose  \cite{Hir}  (1993), Mason and Burdick  \cite{Mas}  (2000); McIsaac and  Ostrowski
 \cite{McI} (2000); Martinez and J. Cortes  \cite{Mar}  (2001);  Trintafyllou et al.  \cite{Tri} (2000);
Morgansen et al.  \cite{Mor} (2001);  Fakuda et al.  \cite{Fak} (2002); Guo et al. \cite{Guo}  (2002); Hawthornee et al.  \cite{Haw}  (2004), and the references therein. 
It was also recognized that sophistication and complexity of design of
bio-mimetic robots give rise to control-theoretic methods, see, e.g., 
 Koiller et al.  \cite{Koi}  (1996); McIsaac and
Ostrowski \cite{McI}  (2000); Martinez and Cortes  \cite{Mar}  (2001);
Trintafyllou et al.  \cite{Tri} (2000); San Martin et al. \cite{San}(2007),
Alouges et al. \cite{Alo} (2008), Sigalotti and  Vivalda  \cite{Sig}  (2009),
and the references therein.  It should be noted however that the above-cited   results  deal with control  problems  in the
framework of ODE's only.

A  number   of attempts were made along these lines to introduce various reduction techniques  to convert swimming model equations into systems of ODE's, namely, by making use of applicable analytical considerations, empiric  observations and experimental data, see, e.g.,
Becker et al  \cite{Bec}  (2003); Kanso et al. \cite{Kan} (2005); San Martin et al.  \cite{San}  (2007); Alogues et al.  \cite{Alo} (2008), and the references therein.

\bigskip
{\bf Problem formulation for 3-$D$ swimming model.}
We consider the following   model, consisting of two {\em coupled} systems of equations: one is a PDE  system--  for the fluid, governed by the 
{\em nonstationary} 3-$D$ Stokes  equation, and the other is an ODE  system-- for the {\em position of the swimming object (or swimmer)} in it:

$$
   \frac{\partial y}{\partial t} \; = \; \nu \Delta y \; + \; F(z, v) \;  -  \;
   \nabla p  \;\;\; \;\;  {\rm in} \;\;\;
   Q_T = \Omega \times (0, T), \;\;  y = (y_1, y_2, y_3)\eqno(1.1)$$
   $$
   {\rm div} \, y = 0 \;\;  {\rm in} \;\;\; Q_T, \;\;\;\; y = 0 \;\;\;\; {\rm in} \;\;
   \Sigma_T = \partial \Omega \times (0, T), \;\;\;\; y \mid_{t = 0} \; = y_0  \;\;{\rm in} \;\Omega,
   $$

$$
\frac{d z_i}{dt} \; = \; \frac{1}{{\rm mes} \, \{ S_i (0) \}}
\mathop{\int}_{S_i (z_i (t))} y (x, t) dx, \;\;\;\; z_i (0) = z_{i0}, \;\;\;\;
i = 1, \ldots, n, \;\;\;\; n > 2,
\eqno(1.2)$$
where for $ t \in [0, T]$:
$$
z (t) = (z_1 (t),   \ldots, z_n (t)), \; z_i (t)  \in R^{3}, \;\; i = 1, \ldots, n, \;\;\;\;  v (t) = (v_1 (t),   \ldots, v_{n-2} (t)) \in R^{n-2},
$$

$$
F(z, v) \; = \;
\sum_{i=2}^{n} [\xi_{i-1}(x, t)k_{i-1}\frac{(\|z_i(t) -
z_{i-1}(t) \|_{R^3} - l_{i-1})}{\|z_i(t) -
z_{i-1}(t)\|_{R^3}}(z_i(t) - z_{i-1}(t))
$$
$$
+ \xi_{i}(x,t)k_{i-1}\frac{(\|z_{i}(t) -
z_{i-1}(t) \|_{{R}^3} - l_{i-1})}{\|z_{i}(t) -
z_{i-1}(t)\|_{{R}^3}}(z_{i-1}(t) -
z_{i}(t))]
$$
$$
+ \sum_{i=2}^{n-1} v_{i-1}(t)\left[\xi_{i-1}(x,t) \left(A_i (z_{i-1}(t)-z_i(t)) \right) (t)  -
\xi_{i+1}(x,t)\frac{\|z_{i-1}(t)-z_i(t)\|^2_{{R}^3}}{\|z_{i+1}(t)-z_i(t)
\|^2_{{R}^3}} \left(B_i (z_{i+1}(t)-z_i(t)\right))(t)  \right]
$$
$$
 - \sum_{i=2}^{n-1}
\xi_{i}(x,t)v_{i-1}(t)\left[\left(A_i  (z_{i-1}(t)-z_{i}(t)\right)) (t) - \frac{
\|z_{i-1}(t)-z_i(t)\|^2_{{R}^3}}{\|z_{i+1}(t)-z_i(t)\|^2_{{R}^3}}\left(B_i  (z_{i+1}(t)-z_{i}(t)\right))(t) \right].
\eqno( 1.3)$$
In the above, $ \Omega  $ is a bounded domain in $ R^3$ with   boundary $ \partial \Omega $ of class $ C^2$, $ \; y (x,t) $  and  $ p (x,t) $ are respectively  the velocity and the pressure of the fluid at point $x = (x_1, x_2, x_3)  \in \Omega $ at time $ t$, while $\nu $ is a kinematic viscosity constant.  
Let us explain the terms in (1.1)-(1.3) in more detail.

{\bf Swimmer:} The swimmer  in (1.1)-(1.3) is modeled as a collection of  $n$  bounded sets $ S_i (z_i(t)), i = 1, \ldots, n$ of non-zero measure (such as balls, parallelepipeds, etc.), identified with the fluid within the space which they occupy. These sets are assumed to be open bounded connected sets symmetric relative to the points $z_i (t)$'s which are their   centers of mass.  The sets  $ S_i(z_i (t) )$'s  are viewed as the given sets $ S_i (0) $'s (``0'' stands of the origin) that have been  shifted to the respective positions $z_i (t)$'s without changing their  orientation in space. Respectively, 
$$
\xi_i (x,t) =
\left\{ \begin{array}{ll}
1, \;\;\;\; &  {\rm if} \;\; x \in S_i (z_i (t)),   \\
0, \;\;\;\; &  {\rm if} \;\;x \in \Omega \backslash S_i (z_i (t)),
\end{array}
\right.   \;\;\;\; i = 1, \ldots, n.
\eqno(1.4)$$
Throughout the paper we assume that each $ S_i(0)$ lies  in a ``small'' neighborhood of the origin of given  radius $r>0$, while $ S_i (a) $ denotes the set $ S_i (0)$ shifted to point $ a$. Denote by 
$$
S^0 = \max_{i = 1, \ldots, n} \{{\rm mes} \, \{ S_i (0) \}\}, \;\;\;\; S_0 = \min_{i = 1, \ldots, n} \{{\rm mes} \, \{ S_i (0) \}\}.
$$

{\bf  Forces:} We assume that these sets are subsequently linked by  forces described by the term$ \; F(z, v) $. No ``actual'' physical links between sets $ S_i (z_i (t))$ are assumed   (i.e., they are assumed to be negligible in terms of affect on the resisting surrounding fluid).
The forces in (1.3) are  {\em internal}, relative to the swimmer -- their sum is zero. We assume that a  force applied to a set $ S (z_i(t))$ acts evenly upon all  its  points, and, as such, it creates an external force  on the  fluid surrounding  $ S (z_i(t))$. 

The structural integrity of the  swimming object  is preserved by the elastic  forces acting  according to Hooke's Law. They act   along the lines connecting the respective adjacent centers $ z_i(t)$'s when the distances between any two adjacent points $ z_{i-1} (t) $ and $ z_{i} (t), i = 2, \ldots,  n$  deviate from the respective given  values
$
 l_{i-1} \; > 0, 
i = 2, \ldots,  n$
as described in the first sum in (1.3).
The  parameters $ \; k_i > 0, \; i = 1, \ldots, n-1 $
characterize the rigidity of the links $ z_{i-1} (t) z_i (t)$, $ i = 2, \ldots, n $. 
The matching  pairs of these forces between $ z_{i-1} (t) $ and $ z_{i} (t) $, and between $ z_{i} (t) $ and $ z_{i+1} (t) $ are shown  on Fig. 1. 

The 2nd  sum in  (1.3) describes the rotation forces about  any on the points $ z_i (t), i = 2, \ldots, n-1$ which make
the adjacent points to rotate about it perpendicular  to the lines connecting the respective  $ z_i(t)$'s. To satisfy the 3rd Newton's Law, these forces  lie in the same plane along with the matching counter-force given in the 3rd sum in (1.3). Respectively, 
$$ 
A_i  = A_i ( z_{i-1}(t), z_i(t)), z_{i+1}(t)  ): \; R^3 \; \rightarrow \; {\rm span} \,\{z_{i - 1}(t)-z_i(t),  
z_{i + 1}(t)-z_i(t) )\}, 
$$
$$
B_i  = B_i ( z_{i-1}(t), z_i(t)), z_{i+1}(t)  ): \; R^3 \; \rightarrow \; {\rm span} \,\{z_{i - 1}(t)-z_i(t), z_{i + 1}(t)-z_i(t)  \}
$$
 denote a {\em nonlinear}  mappings, {\em defined at each moment of time by three vectors $  z_{i-1}(t), z_i(t)), z_{i+1}(t)$},  such that 
for $\;i = 2, \ldots, n-1$:
\begin{itemize}
\item
$$
(z_{i - 1}(t)-z_i(t))^\prime [(A_i \left(z_{i -1}(t)-z_i(t)\right) ) (t)] = 0, \;\;
(z_{i + 1}(t)-z_i(t))^\prime [(B_i \left(z_{i + 1}(t)-z_i(t)\right)) (t)] = 0;
$$
\item
$$
\parallel (A_i \left(z_{i - 1}(t)-z_i(t)\right)) (t) \parallel_{R^3} \; = \; \parallel z_{i - 1}(t)-z_i(t)\parallel_{R^3}, 
$$
$$
\parallel (B_i \left(z_{i + 1}(t)-z_i(t)\right) ) (t) \parallel_{R^3} \; = \; \parallel z_{i + 1}(t)-z_i(t)  \parallel_{R^3};
$$
\item
and the directions of vectors $ (A_i \left(z_{i -1}(t)-z_i(t)\right))(t) $ and $ ( B_i \left(z_{i + 1}(t)-z_i(t)\right))(t)$ are such that they correspond to either folding or unfolding motion of lines $ z_{i-1}(t) z_i (t) $ and $ z_{i}(t) z_{i + 1} (t) $ relative to the point $ z_i(t)$.
\end{itemize}
\par\noindent
The magnitudes
and directions of the rotation forces are determined by the given
coefficients $v_i(t)$, $i = 1, \dots, n - 2$.
The choice of fractional coefficients at terms $A_i (z_{i+1}(t)-z_i(t))$  in (1.3)  ensures that the momentum of swimmer's internal forces is conserved at any $ t  \in (0, T)$ (see calculations in \cite{Kh1} in the 2-$D$  case).
A matching pair of rotation forces, generated by $ z_2 (t)$ for the adjacent points, is shown  on Fig. 1.

\bigskip
{\bf  Swimmer's motion.}  Dynamics of  points $ z_i (t) \xi_i (x,t) , i = 1, \ldots, n$ are determined
by the average motion of the fluid within their respective
supports $ \; S_i (z_i (t))$'s  as described in (1.2).

\bigskip
{\bf  Local and global approach to solutions of  (1.1)-(1.3).} Note that, when the adjacent points in the
swimmer's body share the same position in space, the forcing term $ F $ in (1.3)
and hence  model (1.1)-(1.3) become undefined. While such
situation  {\em mathematically}  seems possible, it does not have  to happen. 
First of all, one can address the issue of local existence of solutions to (1.1)-(1.3) on some 
``small'' time-interval $ (0, T)$, assuming that initially model (1.1)-(1.3) is well-defined in the above sense. This is the primary subject of this paper  (see the next section). 
Then the question of global existence can be viewed as   the issue of suitable  selection of coefficients  $ v_i$'s with the purpose
to ensure that  the aforementioned ill-posed situation is avoided.

In model (1.1)-(1.3)  we chose the fluid governed by the {\em nonstationary} Stokes equation which, along with  its  stationary version, is a typical choice of fluid for  micro-swimmers (the case of Low Reynolds numbers). The empiric reasoning
behind this is that, due to the small size  of swimmer, the inertia terms in the Navier-Stokes equation, containing the 1-st order derivatives in $t$ and $x$, can be omitted, provided that the frequency parameter of  the swimmer at hand  is a quantity of order unity.
However, it was noted that a microswimmer  (e.g., a nano-size robot) may use a rather high frequency of motion, which may justify at least in some cases the need for the term $ y_t$ in  Stokes model equations. In general, it seems reasonable to suggest that the presence of this term (in a number of cases) can provide a better approximation of the Navier-Stokes equation than the lack of it. We also point out that in \cite{Fau1}, \cite{Fau2}, \cite{Pes2}  the full-size Navier-Stokes equation is used for  micro-swimmers.
It also seems  that the methods we use  for the nonstationary Stokes equation\index{nonstationary Stokes equation} (as opposed to stationary Stokes equation), may serve as a natural step  toward the swimming models, based on the Navier-Stokes equation.

\bigskip
{\bf 2. Main result:  Local existence and uniqueness.} Let
$\dot{J}(\Omega)$ denote the set of infinitely differentiable
vector functions  with values in ${R}^3$ which have compact
support in $\Omega$ and are divergence-free, i.e., $\textrm{div}
\phi = 0$ in $\Omega$. Denote by $ J_o(\Omega)$ the closure of
this set in the $ (L^2 (\Omega))^3$-norm and by $ G (\Omega)$
denote the orthogonal complement of $ J_o(\Omega)$ in  $ (L^2
(\Omega))^3$ (see, e.g., \cite{Lad}, \cite{Tem}).
In $\dot{J}(\Omega)$ introduce the scalar product
$$
[\phi_1, \phi_2] \; = \; \mathop{\int}_\Omega \sum_{j = 1}^3 \sum_{i = 1}^3 \phi_{1j_{ x_i}} \phi_{2j_{ x_i}} dx, \;\;\;\;
\phi_1 (x) = ( \phi_{11}, \phi_{12}, \phi_{13}), \;\; \phi_2 (x) = ( \phi_{21}, \phi_{22}, \phi_{23}).
$$
Denote by  $H(\Omega)$ the Hilbert space which is the completion
of $\dot{J}(\Omega)$ in the norm
$$
\|\phi_1 \|_{H(\Omega)} = \sqrt{ \mathop{\int}_\Omega \sum_{j =
1}^3 \sum_{i = 1}^3 \phi^2_{1j_{ x_i}}  dx }.
$$
Everywhere below we will assume the following two assumptions:

\bigskip

{\bf Assumption 2.1.}  {\it For the given $ r>0$, defining the size of sets $ S_i (0) $ in (1,4), assume that
$$
l_{i-1} > 2r, \;\; i = 2,\dots,n; \;\; \overline{S_i}(z_i(0))
\subset \Omega, \;\; \|z_{i,0}-z_{j,0}\|_{{R}^3}
> 2r, \ \ i, j = 1, \dots, n,
\ i \not= j;
\eqno(2.1)$$
and the sets $S_i (0), i = 1, \ldots, n $ are such that
$$
\int_{(S_i (0) \cup S_i (h)) \setminus (S_i(0) \cap S_i (h))}dx
= \int_{\Omega}|\xi_i(x) - \xi_i(x - h)|dx \leq C\|h\|_{{R}^3} \
\ \ \forall h \in B_{h_0}(0)
\eqno(2.2)$$
for some
positive constants $h_0$ and $C$, where $\xi_i(x)$ is the
characteristic function of $S_i(0)$ and $B_{h_0}(0) = \{x \mid \parallel x \parallel_{R^3} < h_0\} \subset {R}^3$}.

Conditions (2.1) mean that at time $t = 0$, any two
sets  $S_i (z_i(0))$ do not overlap, and that the swimmer
lies in  $\Omega$. Condition (2.2) is  a
regularity assumption of Lipschitz type regarding the shift of the
set $S(0)$. It is satisfied, for instance, for balls and
parallelepipeds.

{\bf Assumption 2.2.}  {\it Assume that within some $ (R^3)^n$ neighborhood $ {\mathbf G} (z(0)) \subset (R^3)^n$ of the initial datum in (1.2) the mappings $ A_i$ and $ B_i $ are Lipschitz for all $ i = 2, \ldots, n-1$ in the following sense:
$$
\parallel A_i (a_{i-1}, a_i, a_{i+1}) \left(a_{i - 1}-a_i\right) -  A_i (b_{i-1}, b_i, b_{i+1}) \left(b_{i - 1}-b_i\right)) \parallel_{R^3} \; 
$$
$$
\leq  \; L \left\{
\parallel  a_{i - 1}-b_{i-1} ) \parallel_{R^3}  +  \parallel   a_{i }-b_{i}  \parallel_{R^3}  + \parallel   a_{i +1}-b_{i+1}  \parallel_{R^3} \right\}, 
$$
$$
\parallel B_i (a_{i-1}, a_i, a_{i+1})\left(a_{i + 1}-a_i\right) -  B_i (b_{i-1}, b_i, b_{i+1})\left(b_{i + 1}-b_i\right) \parallel_{R^3} \; 
$$
$$
\leq  \; L \left\{
\parallel  a_{i - 1}-b_{i-1} ) \parallel_{R^3}  +  \parallel   a_{i }-b_{i}  \parallel_{R^3}  + \parallel   a_{i +1}-b_{i+1}  \parallel_{R^3} \right\}, 
$$
for any $ a_{i \pm 1}, a_i, b_{i \pm 1}, b_i \in  {\mathbf G} (z(0) $, where $ L > 0$ is a constant and $ A_i$'s and $B_i$'s are defined as in the introduction by three respective vectors $a_{i-1}, a_i, a_{i+1}$.
}

Assumption 2.2 can be satisfied if, e.g., the points $ a_{i-1} = z_{i-1} (0), a_i = z_{i} (0), a_{i+1} = z_{i+1}(0), i = 2, \ldots, n-1$ are not on the same line and the mappings 
$ A_i $ and $ B_i$ are selected by making use of Gramm-Schmidt orthogonalization procedure for vectors $ a_{i-1} - a_i$ and $  a_{i+1} -a_i$.
Alternatively, we can define $ A_i$'s and $B_i$'s, making use of the cross-product:
$$
(A_i (z_{i-1} (t) - z_i (t))(t) \; = \;  \; e_1 (t) \parallel z_{i-1} (t) - z_i (t)) \parallel_{R^3}, 
$$
$$
(B_i (z_{i+1} (t) - z_i (t))(t) \; = \;  \; e_2(t)  \parallel z_{i+1} (t) - z_i (t)) \parallel_{R^3}, \;\;\;\;
e_i (t) = \frac{v_i(t) }{\parallel v_i (t) \parallel_{R^3}}, \;\; i = 1, 2,
$$
$$
v_1 (t) = (z_{i-1} (t) - z_i (t)) \times  [(z_{i-1} (t) - z_i (t)) \times (z_{i+1} (t) - z_i (t))], 
$$
$$
v_2 (t) =   [(z_{i+1} (t) - z_i (t))  (z_{i+1} (t) - z_i (t))] \times (z_{i+1} (t) - z_i (t)) .
$$

\bigskip

Here is the main result of this paper.

{\bf Theorem 2.1}. \  {\it Let $y_0 \in H(\Omega) $; $T
> 0$; $k_i > 0$, $i=1,\dots,n-1$; $v_i \in L^\infty(0, T)$, $i = 1, \dots, n - 2$;
and $z_{i} (0) \in \Omega$, $i = 1, \dots, n$, and let Assumptions
2.1 and 2.2  hold.  Then there exists a $T^* = T^*(z_{1}(0), \dots, z_{n}(0), \|v_1\|_{L^\infty(0, T)}$, $\dots, \|v_{n - 2}\|_{L^\infty(0,
T)}, \Omega) \in (0, T)$ such that system (1.1) - (1.3) admits a unique solution $\left\{y, p, z\right\}$ on $(0,
T^*)$, $\left\{y, \nabla p, z\right\} \in L^2 (0, T^*; J_o(\Omega)) \times
L^2 (0, T^*; G(\Omega)) \times [C([0, T^*]; {R}^3)]^n$.  Moreover, $y \in
C([0, T^*]; H(\Omega))$, $y_t, y_{x_i x_j}, p_{x_i} \in
(L^2(Q_{T^*}))^3$, where $i, j = 1, 2, 3$, and equations (1.1) and
(1.2) are satisfied almost everywhere, while Assumptions 2.1 and 2.2 hold in $ [0, T_*]$.
}

\bigskip
{\bf Remark 2.1.}
\begin{itemize}{}
\item
The fact that  conditions Assumptions 2.1 and 2.2 hold in $ [0, T_*]$  implies  that  we are able to guarantee that within $[0, T^*]$
no parts of the swimmer's body will ``collide'',
and simultaneously, that it stays strictly inside of $\Omega$. These
conditions allow us to maintain the mathematical and physical
wellposedness of model (1.1) - (1.3).
\item
As it will follow from the proof below,  
Theorem $2.1$  allows
further extension of the solutions to (1.1) - (1.3) in time as long as Assumptions 2.1 and 2.2 continue to hold.  This
depends on the choice of parameters  $v_1(t), \dots, v_{n-2}(t)$.
\end{itemize}

Our plan to prove Theorem $2.1$ is to proceed stepwise as follows:

\begin{itemize}
\item In Section $3$ we discuss the existence and uniqueness of the solutions to
the decoupled version of (1.2).  
\item
In Section 4 we will introduce three
continuous mappings for the decoupled version of the system (1.1) -
(1.3).
\item
 In Section 5 we will apply a fixed point argument to prove
Theorem 2.1. 
\end{itemize}

In the proofs below we employ the methods introduced in \cite{Kh1} to investigate the wellposedness of  the 2-$D$ version of model (1.1)-(1.3), 
modifying and extending them  to the 3-$D$ case.

Without loss of generality, we will further assume that system (1.1)-(1.3) and all respective auxiliary systems below are considered on the time-intervals whose lengths are smaller than 1. 

\bigskip
{\bf 3. Preliminary results.} 
Introduce the following decoupled version of system (1.2):
$$
\frac{dw_i}{dt} =
\frac{1}{\textrm{mes}(S_i(0))}\int_{S_i(w_i(t))}u(x, t)dx, \ \ \
w_i(0) = z_{i, 0}, \ \ \ i = 1, \dots, n,
\eqno(3.1)$$
where $u(x, t)$ is some given function. Denote $ w (t) = (w_1(t), \ldots, w_n (t)) $.

{\bf Lemma $3.1$}. \ \  {\it Let $T > 0$ and $u \in (L^2(0,T;L^\infty(\Omega)))^3$ be given.
Then there is a $T^* \in (0,T)$ such that system (3.1) has a unique
solution in $C([0, T^*]; {R}^3)$ satisfying Assumptions 2.1 and 2.2 with $ w (t) $ in place of $ z(t)$, if they hold at time $ t = 0$.}

{\em Proof}. We will use the contraction principle
to prove existence and uniqueness.
Below the values of  $h_0$, $C$ are taken from (2.2).

Select $ T_0$ to satisfy the following inequalities: 
$$
0 < T_0 < \min\left\{\frac{\textrm{mes}(S_0)h_0^2}{4\|u\|_{(L^2(Q_T))^3}^2},
\frac{(\textrm{mes}(S_0))^2}{C^2\|u\|_{(L^2(0,T;L^\infty(\Omega)))^3}^2} 
T, 1\right\}
$$
$$
\leq  \min_{i = 1, \ldots, n} \min\left\{\frac{\textrm{mes}(S_i(0))h_0^2}{4\|u\|_{(L^2(Q_T))^3}^2},
\frac{(\textrm{mes}(S_i(0)))^2}{C^2\|u\|_{(L^2(0,T;L^\infty(\Omega)))^3}^2} 
T, 1\right\}.
\eqno(3.2)$$
Let for any given $ p \in C([0,T_0];{R}^3)$: 
$$
B_{h_0/2}(p) = \left\{z \in C([0,T_0];{R}^3) \mid
\|z - p\|_{C([0,T_0];{R}^3)} \leq \frac{h_0}{2}\right\} \subset
C([0,T_0];{R}^3).
$$
For each $i=1,\dots,n$, define a mapping
$D_i: B_{h_0/2}(z_{i,0}) \longrightarrow C([0,T_0];{R}^3)$ by
$$
D_i(w_i(t)) =
z_{i,0}+\frac{1}{\textrm{mes}(S_i(0))}\int_{0}^{t}\int_{S_i(w_i(\tau))}u(x,\tau)dx d\tau.
$$
Then we can derive that:
$$
\|D_i(w_i(t))\|_{{R}^3} 
\leq \ \|z_{i,0}\|_{{R}^3} +
\frac{\sqrt{T_0}}{\sqrt{\textrm{mes}(S_i (0))}}\|u\|_{(L^2(Q_T))^3} \
\forall \ t \in [0, T_0], \ i = 1, \dots, n.
\eqno(3.3{\rm a})$$
Similarly, in view of (3.2):
$$
\|D_i(w_i(t))-z_{i,0}\|_{C([0,T_0];{R}^3)} \leq
\frac{\sqrt{T_0}}{\sqrt{\textrm{mes}(S_i (0))}}\|u\|_{(L^2(Q_T))^3} <
\frac{h_0}{2}.
\eqno(3.3b)$$
Thus, $D_i$ maps $B_{h_0/2}(z_{i,0})$ into itself for each $i=1,\dots,n$, where $z_{i,0}$ is treated as a constant function.

Let $w_i^{(1)}(t),w_i^{(2)}(t) \in B_{h_0/2}(z_{i,0})$ and $\xi(x,S)$ denote the
characteristic function of a set $S \subset {R}^3$.  
Then, making use of  (2.2), we obtain:
$$
\|D_i(w_i^{(1)}(t)) - D_i(w_{i}^{(2)}(t))\|_{{R}^3}
=
$$
$$
\frac{1}{\textrm{mes}(S_i(0))}\left\|\int_{0}^{t}\int_{S_i(w_i^{(1)}(t))}u(x,
\tau)dxd\tau - \int_{0}^{t}\int_{S_i(w_{i}^{(2)}(\tau))}u(x, \tau)dxd\tau\right\|_{{R}^3}
$$
$$
= \frac{1}{\textrm{mes}(S_i(0))}\left\|\int_{0}^{t}\int_{\Omega}u(x,
\tau)\left(\xi(x, S_i(w_i^{(1)}(t))) - \xi(x,
S_i(w_{i}^{(2)}(\tau)))\right)dxd\tau \right\|_{{R}^3}
 $$
$$
\leq
\frac{C\sqrt{T_0}}{\textrm{mes}(S_i(0))}\|u\|_{(L^2(0,T;L^\infty(\Omega)))^3}
\|w_i^{(1)}(t)-w_i^{(2)}(t)\|_{C([0,T_0]; {R}^3)}, \ \forall \
t \in [0, T_0], \ i = 1, \dots, n.
\eqno(3.4)$$
Therefore, after maximizing the left-hand side of (3.4) over
$[0,T_0]$, we conclude:
$$
\|D_i(w_i^{(1)}(t)) -
D_i(w_{i}^{(2)}(t))\|_{C([0,T_0];{R}^3)}\leq
\frac{C\sqrt{T_0}}{\textrm{mes}(S_i(0))}\|u\|_{(L^2(0,T;L^\infty(\Omega)))^3}
\|w_i^{(1)}-w_i^{(2)}\|_{C([0,T_0];{R}^3)}.
\eqno(3.5)$$
Hence, in view of (3.2),
$$
\frac{C\sqrt{T_0}}{\textrm{mes}(S_i(0))}\|u\|_{(L^2(0,T;L^\infty(\Omega)))^3} < 1,
$$
it follows from (3.5) that $D_i$ is a contraction mapping on
$B_{h_0/2}(z_{i,0})$ for each $i=1,\dots,n$. Therefore, there exist
unique $w_i(t) \in C([0,T_0];{R}^3), i =  1, \ldots, n$, such that
$D_i(w_i(t))=w_i(t)$, i.e.,
$$
w_i(t) =
z_{i,0}+\frac{1}{\textrm{mes}(S_i(0))}\int_{0}^{t}\int_{S_i(w_i(\tau))}u(x,\tau)dxd\tau,
\ i = 1,\dots,n,
\eqno(3.6)$$
which yields (3.1).

{\bf On restrictions  (2.3):}  Estimates  (3.3a-b)  imply that we may select a $T^*
\in (0,T_0)$ such that for any $t \in [0, T^*]$, all $w_i(t)$'s will
stay ``close enough'' to their initial values $z_{i,0}$'s to guarantee 
that Assumptions 2.1 and 2.2  holds for $w_i(t)$, $i=1,\dots,n$. This ends the proof of
Lemma 3.1.

\bigskip

{\bf 4. Decoupled solution mappings.}
Let $ {\cal B}_q(0)$ denote a closed ball of radius $q$ (its value  will be
selected in Section $5$) with center at the origin in the Banach space
$ L^2 (0, T; J_o(\Omega)) \bigcap L^2(0, T; (H^2(\Omega))^3)$ endowed with the norm
of $L^2(0, T; (H^2(\Omega))^3)$:
$$
{\cal B}_q(0) = \left\{\phi \in L^2 (0, T; J_o(\Omega)) \bigcap
L^2(0, T; (H^2(\Omega))^3) \ | \ \|\phi\|_{L^2(0, T;
(H^2(\Omega))^3)} \leq q\right\},
$$
where $ H^2 (\Omega) = \{\phi \mid \phi, \phi_{x_i}, \phi_{x_i x_j} \in L^2 (\Omega), \, i,j = 1,2,3\}$.

Note that
$H^2(\Omega)$ is continuously embedded into $C(\bar{\Omega})$,
and thus $L^2(0, T; (H^2(\Omega))^3)$ is continuously embedded into
$(L^2(0, T; L^\infty(\Omega)))^3$. This yields the estimate
$$
\|\phi\|_{(L^2(0, T; L^\infty(\Omega)))^3} \leq
K\|\phi\|_{L^2(0, T; (H^2(\Omega))^3)} \ \ \ \textrm{for some} \ K > 0.
\eqno(4.1)$$
This implies that Lemma $3.1$ holds for
any $u \in {\cal B}_q(0)$.

\bigskip
{\bf 4.1. Solution mapping for $ z_i(t), i = 1, \dots, n$.}
We now intend to show that the operator
$$ {\mathbf A}: {\cal B}_q(0) \longrightarrow [C([0, T]; {R}^3)]^n, \ \ \ \ \
{\mathbf A}u = w = (w_1, \dots, w_n),$$ where the $w_i$'s solve (3.1), is
continuous and compact if $T > 0$ is sufficiently small.

\bigskip
{\bf Continuity.} \ Let $u^{(1)}, u^{(2)} \in {\cal B}_q(0)$ with $T_1$ in
place of $T$, where $T_1 > 0$ satisfies assumptions in the proof of Lemma 3.1 with $T_1$ in place of
$T^*$. Define ${\mathbf A} u^{(j)} = w^{(j)} = (w_1^{(j)}, \dots, w_n^{(j)})$
for $j = 1, 2$. To show $ {\mathbf A} $ is continuous, we will evaluate
$$
\|{\mathbf A} u^{(1)} - {\mathbf A} u^{(2)}\|_{[C([0, T]; {R}^3)]^n}\
$$
term-by-term.  To this end, similar to (3.4), we have the following
estimate:
$$
\|w_{i}^{(1)}(t) - w_{i}^{(2)}(t)\|_{{R}^3}
$$
$$
=\left\|\frac{1}{\textrm{mes}(S_i(0))}\int_{0}^{t}\int_{S(w_i^{(1)}(\tau))}u^{(1)}(x,\tau)
dxd\tau -
\frac{1}{\textrm{mes}(S_i(0))}\int_{0}^{t}\int_{S(w_i^{(2)}(\tau))}u^{(2)}(x,\tau)
dxd\tau\right\|_{{R}^3}
$$
$$
= \frac{1}{\textrm{mes}(S_i(0))} \parallel \int_{0}^{t}\int_{S(w_i^{(1)}(\tau))}u^{(1)}(x,\tau)
dxd\tau - \int_{0}^{t}\int_{S(w_i^{(1)}(\tau))}u^{(2)}(x,\tau)
dxd\tau
$$
$$
 +
\int_{0}^{t}\int_{S(w_i^{(1)}(\tau))}u^{(2)}(x,\tau) dxd\tau -
\int_{0}^{t}\int_{S(w_i^{(2)}(\tau))}u^{(2)}(x,\tau)
dxd\tau \parallel_{{R}^3}
$$
$$
\leq \frac{1}{\textrm{mes}(S_i(0))} (\left\|\int_{0}^{t}\int_{\Omega}(u^{(1)}(x,\tau)
- u^{(2)}(x,\tau))\xi(x,
S(w_i^{(1)}(\tau)))dxd\tau\right\|_{{R}^3}
$$
$$
 +
\left\|\int_{0}^{t}\int_{\Omega}u^{(2)}(x,\tau)(\xi(x,
S(w_i^{(1)}(\tau)))
- \xi(x, S(w_i^{(2)}(\tau))))dxd\tau\right\|_{{R}^3} )
$$
$$
\leq \frac{\sqrt{T_1}}{\sqrt{\textrm{mes}(S_i(0))}}\|u^{(1)} -
u^{(2)}\|_{(L^2(Q_{T_1}))^3} +
\frac{C\sqrt{T_1}}{\textrm{mes}(S_i(0))}\|u^{(2)}\|_{(L^2(0, T_1;
L^\infty(\Omega)))^3}\|w_i^{(1)} - w_i^{(2)}\|_{C([0, T_1]; {R}^3}
$$
$$
\leq \frac{\sqrt{T_1}}{\sqrt{\textrm{mes}(S_0)}}\|u^{(1)} -
u^{(2)}\|_{(L^2(Q_{T_1}))^3} +
\frac{C\sqrt{T_1}}{\textrm{mes}(S_0)}\|u^{(2)}\|_{(L^2(0, T_1;
L^\infty(\Omega)))^3}\|w_i^{(1)} - w_i^{(2)}\|_{C([0, T_1]; {R}^3}.
\eqno(4.2)$$
Recall from (4.1) that
$$
\|u^{(2)}\|_{(L^2(0, T_1; L^\infty(\Omega)))^3} \leq Kq.
$$
So,
select a $T > 0$ as follows:
$$
0 < T <
\textrm{min}\left\{\left(\frac{\textrm{mes}(S_0)}{CKq}\right)^2,
T_1\right\}.
\eqno(4.3)$$
Hence, replacing $T_1$ in (4.2)
with $T$ satisfying (4.3) and maximizing the left-hand side of (4.2)
over $[0, T]$, we obtain:
$$
\|w_i^{(1)} - w_i^{(2)}\|_{C([0, T]; {R}^2)}  \leq
\frac{\sqrt{T}}{\sqrt{\textrm{mes}(S_0)}}\|u^{(1)} -
u^{(2)}\|_{(L^2(Q_{T}))^3}
$$
$$
+ \frac{CKq\sqrt{T}}{\textrm{mes}(S_0)}\|w_i^{(1)} -
w_i^{(2)}\|_{C([0, T]; {R}^3)}.
$$
In view of (4.3), if $w_i^{(1)}(t) \not= w_i^{(2)}(t)$ on $[0, T]$,
then the above implies:
$$
0 < \left(1 -
\frac{CKq\sqrt{T}}{\textrm{mes}(S_0)}\right)\|w_i^{(1)} -
w_i^{(2)}\|_{C([0, T]; {R}^3)} \leq
\frac{\sqrt{T}}{\sqrt{\textrm{mes}(S_0)}}\|u^{(1)} -
u^{(2)}\|_{(L^2(Q_{T}))^3}.
$$
 Thus, it follows that
$$
\|w_i^{(1)} - w_i^{(2)}\|_{C([0, T]; {R}^3)} \leq
\frac{\sqrt{T\textrm{mes}(S_0)}}{\textrm{mes}(S_0) -
CKq\sqrt{T}}\|u^{(1)} - u^{(2)}\|_{(L^2(Q_T))^3}.
\eqno(4.4)$$

Therefore, (4.3) and (4.4) imply that for every $u^{(1)}$, $u^{(2)}
\in {\cal B}_q(0)$,
$$
\|{\mathbf A} u^{(1)} - {\mathbf A} u^{(2)}\|_{[C([0, T]; {R}^3)]^n} \leq
\frac{\sqrt{nT\textrm{mes}(S_0)}}{\textrm{mes}(S_0) -
CKq\sqrt{T}}\|u^{(1)} - u^{(2)}\|_{L^2(0,T;(H^2(\Omega))^3)}.
$$
So, the operator ${\mathbf A} $ is continuous on ${\cal B}_q(0)$ for sufficiently
small $T$ as in (4.3).

{\bf Compactness.} \ Furthermore, to show that ${\mathbf A} $ is compact, we
will show that ${\mathbf A} $ maps any sequence in ${\cal B}_q(0)$ into a sequence in
$[C([0, T];{R}^3)]^n$ which contains a convergent
subsequence. To this end, consider any sequence
$\left\{u^{(j)}\right\}_{j=1}^{\infty}$ in ${\cal B}_q(0)$. Using (3.6)
with $w_i^{(j)}$ and $u^{(j)}$ in place of $w_i$ and $u$, construct  the
sequence $\left\{w_i^{(j)}\right\}_{j=1}^{\infty}$, $i = 1, \dots,
n$.

Let us now show that $\left\{w_i^{(j)}\right\}_{j=1}^{\infty}$ is
uniformly bounded and equicontinuous.  Indeed, applying (4.1) to an
estimate like (3.3a) and then maximizing over $[0, T]$ yields:
$$
\|w_i^{(j)}\|_{C([0, T]; {R}^3)} \leq \max_{i = 1, \dots,
n}\left\{\|z_{i, 0}\|_{{R}^3}\right\} +
\frac{q\sqrt{T}}{\sqrt{\textrm{mes}(S_0)}}.
\eqno(4.5)$$
To show equicontinuity, consider any $t, t+h \in [0, T]$, e.g., when
$h>0$. Then for $i = 1, \dots, n$:
$$
\|w_i^{(j)}(t + h) - w_i^{(j)}(t)\|_{{R}^3} = \
 \ \frac{1}{\textrm{mes}(S_i(0))}\left\|\int_{t}^{t +
h}\int_{S(w_{i}^{(j)}(\tau))}u^{(j)}(x, \tau)dxd\tau\right\|_{{R}^3}
$$
$$
\leq \ \frac{\sqrt{h}}{\sqrt{\textrm{mes}(S_i(0))}}\|u^{(j)}\|_{(L^2(Q_T))^3}
\leq \frac{q\sqrt{h}}{\sqrt{\textrm{mes}(S_0)}} \quad j = 1,\dots,
$$
which implies the equicontinuity of
$\left\{w_i^{(j)}\right\}_{j=1}^{\infty}$, $i = 1, \dots, n$ (the
case $h < 0$ is similar). Therefore, by Ascoli's Theorem,
$\left\{{\mathbf A} u^{(j)}\right\}_{j=1}^{\infty}$ contains a convergent
subsequence in $[C([0, T]; {R}^3)]^n$, i.e., ${\mathbf A} $ is compact
 on $ {\cal B}_q(0)$.

\bigskip
{\bf 4.2. Solution mapping for
decoupled Stokes equations.} Now, consider the following
decoupled Stokes initial boundary value problem:
$$
\frac{\partial y_\ast}{\partial t} - \nu \Delta y_\ast + \nabla p_\ast =
f(x, t) \ \ \ \textrm{in} \ \ Q_T,
$$
$$
\textrm{div} \ y_\ast = 0 \ \ \ \textrm{in} \ \
Q_T, \ \ \ \ y_\ast = 0 \ \ \ \textrm{in} \ \ \Sigma_T, \ \ \
y_\ast|_{t = 0} = y_0 \in H(\Omega).
\eqno(4.6)$$

For any $f \in (L^2(Q_T))^3$, it is known that the boundary value
problem (4.6) possesses a unique solution $y_\ast$ in $ L^2 (0, T; J_o(\Omega))
\bigcap L^2(0, T; (H^2(\Omega))^3)$ with the properties described in
Theorem $2.1$ (see, e.g. [19], [27]).  Moreover, (see, e.g. (7) on
p. 79 and  (50) on p. 65 in [19]), there
is a positive constant $L$ such that:
$$
\|y_*\|_{L^2(0, T; (H^2(\Omega))^3)}^2 \leq L\|y_0\|_{H(\Omega)}^{2}
+ L\int_{Q_T}\|f(\cdot, \tau)\|_{{R}^3}^2d\tau.
\eqno(4.7)$$
Thus it follows that, given $y_0$, the operator
$${\mathbf B}: (L^2(Q_T))^3 \longrightarrow
L^2 (0, T; J_o(\Omega)) \bigcap L^2(0, T; (H^2(\Omega))^3), \quad {\mathbf B}f = y_\ast$$ is continuous.

\bigskip
{\bf 4.3. The  term
${\mathbf F(z, v)}$.} Let $T > 0$ be as in (4.3).  Given $u \in
{\cal B}_q(0)$, let $F_\ast(w), w = (w_1, \ldots, w_n)$ denote the value of $F(z, v)$ in (1.3),
with  $w_i$'s from (3.6) in  $z_i$'s, respectively. Consider the operator
$$
{\mathbf F} : \ [C([0, T]; {R}^3)]^n \longrightarrow (L^2(Q_T))^3, \ \\ \;\; {\mathbf F}w = F_\ast(w).
$$
We will show that ${\mathbf F}$ is
continuous, but first we evaluate $\|F_\ast(w)\|_{(L^2(Q_T))^2}$.
To do that, we will use the standard algebraic  transformation technique  similar to that  used for the 2-$D$ swimming model in \cite{Kh1}, \cite{Kh2}. Therefore, we omit some of the details below.

Let $P(T)$ denote the upper bound in (4.5):
$$
P(T) = \max_{i = 1, \dots, n}\left\{\|z_{i, 0}\|_{{R}^3}\right\} +
\frac{q\sqrt{T}}{\sqrt{\textrm{mes}(S_0)}}.
\eqno(4.8)$$
Then,
using (4.5) with $w_i$ in place of $w_i^{(j)}$, we can evaluate
the first term in square brackets in the first line of (1.3) as
follows:
$$
\left\|\xi_{i}(x,\tau)k_{i-1}\frac{(\|w_{i}(t) -
w_{i-1}(t)\|_{{R}^3} - l_{i-1})}{\|w_{i}(t) -
w_{i-1}(t)\|_{{R}^3}}(w_{i}(t) -
w_{i-1}(t))\right\|_{{R}^3}
$$
$$
\leq \; C_o \left(n, r, P(T), \max_{i = 1, \dots,
n-1}\left\{k_i\right\} , \max_{i = 1, \dots n -
1} \left\{l_i\right\} \right),
$$
where  and below $C_0 = C_0 \left(n, r, P(T), \max_{i = 1, \dots,
n-1}\left\{k_i\right\} , \max_{i = 1, \dots n -
1} \left\{l_i\right\} \right) $ denotes a positive generic constant, continuous in $ r $, $k_i$'s, $l_i$'s,  and $ P(T)$  (it can have a  different concrete value for different expressions below). Along  similar estimates for other terms in (1.3), we obtain:
$$
\|F_\ast(w)\|_{(L^2(Q_T))^3} \; \leq \; C_o \sqrt{T\textrm{mes}(\Omega)}  \left(1 \; + \;  \max_{i = 1, \dots,
n-2}\left\{\|v_{i}\|_{L^{\infty}(0, T)}\right\}  \right)  \quad
\forall \ u\in {\cal B}_q(0).
\eqno(4.9)$$

We will now show that ${\mathbf F}$ is a continuous operator. Let $w^{(1)},
w^{(2)} \in [C([0,T];{R}^3)]^n$.  Consider (1.3) with
$w_i^{(j)}$ in place of $z_i$, $i = 1, \dots, n$, $j=1,2$.

Making use of estimate like (3.3b), without loss of generality we can assume that $ T$ is small enough to ensure that Assumptions 2.1 and  2.2 holds (as stated in Lemma 3.1).
Then, once again, similar to the 2-$D$ case model in \cite{Kh1}, we can  evaluate $\|F_*(w^{(1)})-F_*(w^{(2)})\|_{{R}^3}$
along the standard algebraic transformations and making use of Assumption 2.2, which will result in the following formula:
$$
\|F_*(w^{(1)}) - F_*(w^{(2)})\|_{(L^2(Q_T))^3} \leq
M(T)\sqrt{T\textrm{mes}(\Omega)}\|w^{(1)} - w^{(2)}\|_{[C([0, T];{R}^3)]^n},
\eqno(4.10)$$
where
$$
M(T) = M(T, n, r, P(T), \max_{i = 1, \dots,
n-1}\left\{k_i\right\} , \max_{i = 1, \dots n -
1} \left\{l_i\right\} ) \; \left(1 \; + \;  \max_{j=1,\dots,n-2}\left\{\|v_j\|_{L^\infty(0, T)}\right\}
\right),
\eqno(4.11)$$
is defined similar to $ C_o$. Hence ${\mathbf F}$ is a continuous operator for $T$ satisfying (4.3).

\bigskip

{\bf 5. Proof of Theorem 2.1.}
Before proceeding to the proof of Theorem 2.1, we summarize the
main results of the previous section.  In Section 4, we proved
that for sufficiently small $T > 0$ satisfying (4.3), the operators
$$
{\mathbf A} : {\cal B}_q(0) \longrightarrow [C([0, T];{R}^3)]^n, \quad {\mathbf A} u = w =
(w_1, \dots, w_n),
$$
$$
{\mathbf F}: [C([0, T]; {R}^3)]^n \longrightarrow (L^2(Q_T))^3, \quad {\mathbf F}w =
F_\ast(w),
$$
and
$$
{\mathbf B}: (L^2(Q_T))^3 \longrightarrow L^3 (0, T; J_o(\Omega))\bigcap
L^2(0, T; (H^2(\Omega))^3), \quad {\mathbf B}f = y_\ast
$$
are all continuous.  ${\mathbf A}$ is also compact. As a result,
$$
 {\mathbf B}{\mathbf F}{\mathbf A} : \ B_q(0) \longrightarrow L^2 (0, T; J_o(\Omega)) \bigcap L^2(0, T; (H^2(\Omega))^3),
\quad {\mathbf B}{\mathbf F}{\mathbf A}  u = y_*
$$
is continuous and compact.

\bigskip
{\bf 5.1. Existence:  Fixed point argument.}
Select the value of $q$ to be any positive number larger than
$\sqrt{L}\|y_0\|_{H(\Omega)}$, and choose
$T > 0$ as in (4.3) so that Lemma $3.1$ holds.  In view of (4.7)
(with $F_\ast$ in place of $f$), we select $T^* \in (0,
\min\left\{T, 1\right\})$ small enough so that the continuous and
compact operator ${\mathbf B} {\mathbf F}{\mathbf A} $ maps the closed ball $B_q(0)$ into itself.

Re-wrute  (4.9) as follows:
$$
\|F_*(w)\|_{(L^2(Q_T))^3} < C_1\sqrt{T},
$$
 where the positive constant $C_1$ is independent of $T$.
Select $T^*$ so that:
$$
0 < T^* < \textrm{min}\left\{\frac{q^2 -
L\|y_0\|_{H(\Omega)}^2}{LC_1^2}, T, 1\right\}.
\eqno(5.1)$$
 Since $T^*$ satisfies (4.3), if we replace $T$ with
$T^*$ in (4.7), (4.9) - (4.11), we obtain:
$$
\|{\mathbf B}{\mathbf F}{\mathbf A} u\|_{L^2(0, T^*; (H^2(\Omega))^3)}^2 \leq
L\|y_0\|_{H(\Omega)}^2 +
L\|F_\ast(w)\|_{(L^2(Q_{T^*}))^3}^2
$$
$$
< L\|y_0\|_{H(\Omega)}^2 + LC_1^2T^* < L\|y_0\|_{H(\Omega)}^2 + q^2
- L\|y_0\|_{H(\Omega)}^2 = q^2.
$$
Hence, ${\mathbf B}{\mathbf F}{\mathbf A} $ maps
${\cal B}_q(0)$ into itself if (5.1) is satisfied.

Thus, by Schauder's Fixed Point Theorem, ${\mathbf B} {\mathbf F}{\mathbf A} $ has a fixed point $y$
which is a solution of the system (1.1) - (1.3), and which satisfies
all of the requirements of Theorem $2.1$. As usual, we may select
$\nabla p$ in $ L^2 (0, T^*; G(\Omega))$ to complement the solution $y \in
L^2 (0, T^*; J_o(\Omega)) $ in Theorem $2.1$. This completes the proof of
existence for Theorem 2.1.

\bigskip
{\bf 5.2. Uniqueness.} To prove that the solution found in Section 5.1
is unique, we will argue by contradiction.  Namely, suppose, e.g.,
that there are two different solutions
$$
\left\{z^{(1)} = (z_1^{(1)},
\dots, z_n^{(1)}), y^{(1)}, p^{(1)}\right\}$$
 and
$$
\left\{z^{(2)} =
(z_1^{(2)}, \dots, z_n^{(2)}), y^{(2)}, p^{(2)}\right\}
$$
 to (1.1)-(1.3), satisfying the properties described in Theorem 2.1 on some
time interval $[0, T]$, where $T$ satisfies inequality (5.1).
Without loss of generality, we assume the two solutions are
different right from $t = 0$.

By (4.4), with $z_i^{(j)}$ in place of $w_i^{(j)}$, $i = 1, \dots,
n$, and $y^{(j)}$ in place of $u^{(j)}$, both for $j = 1, 2$, we see
that for any $T_0 \in (0, T]$, and each $i=1,\dots,n$:
$$
\|z_i^{(1)} - z_i^{(2)}\|_{C([0, T_0]; {R}^3)} \leq
\frac{\sqrt{T_0\textrm{mes}(S_0)}}{\textrm{mes}(S_0) -
CKq\sqrt{T_0}}\|y^{(1)} - y^{(2)}\|_{(L^2(Q_{T_0}))^3}.
\eqno(5.2)$$

Let us now evaluate $\|y^{(1)} - y^{(2)}\|_{(L^2(Q_{T_0}))^3}$. Note
that $(y^{(1)} - y^{(2)})$ satisfies the following Stokes
initial-value problem:
$$
\frac{\partial(y^{(1)} - y^{(2)})}{\partial t} = \nu \Delta(y^{(1)} - y^{(2)}) +
(F(z^{(1)}, v) - F(z^{(2)}, v)) - \nabla(p^{(1)} - p^{(2)}) \ \ \
\textrm{in} \ Q_{T_0},
$$
$$
\textrm{div} \ (y^{(1)} - y^{(2)}) = 0 \ \ \ \textrm{in} \ Q_{T_0}, \ \ \ \
(y^{(1)} - y^{(2)}) = 0 \ \ \ \textrm{in} \ \Sigma_{T_0}, \ \ \
(y^{(1)} - y^{(2)})|_{t = 0} = 0.
$$

According to (4.7) we have:
$$
\|y^{(1)} - y^{(2)}\|_{(L^2(Q_{T_0}))^3}^2 \leq
L\int_{0}^{T_0}\int_{\Omega}\|F(z^{(1)}, v) - F(z^{(2)},
v)\|_{{R}^3}^{2}dxdt.
\eqno(5.3)$$
In turn, similar to (4.10):
$$
\|F(z^{(1)}, v) - F(z^{(2)}, v)\|_{{R}^3} \leq N(T)\sum_{j =
1}^{n}\|z_j^{(1)} - z_{j}^{(2)}\|_{C([0, T_0]; {R}^3)}.
\eqno(5.4)$$
for some  $N(T)$ is nonincreasing at $T \rightarrow 0^+$.  Hence,
combining (5.2) - (5.4) yields:
$$
\|y^{(1)} - y^{(2)}\|_{(L^2(Q_{T_0}))^3} \leq
\frac{nN(T)T_0\sqrt{L\textrm{mes}(S_0)\textrm{mes}(\Omega)}}{\textrm{mes}(S(0))
- CKq\sqrt{T_0}}\|y^{(1)} - y^{(2)}\|_{(L^2(Q_{T_0}))^3}.
\eqno(5.5)$$
Now, select $T_0$ as follows:
$$
0 < T_0 < \min\left\{\frac{\textrm{mes}(S_0)}{4nN(T)\sqrt{L\textrm{mes}(S_0)\textrm{mes}(\Omega)}},
\ \left(\frac{\textrm{mes}(S_0)}{2CKq}\right)^2, \
T\right\}.
\eqno(5.6)$$
This choice of $T_0$ implies that the following inequality holds:
$$
nN(T)T_0\sqrt{L\textrm{mes}(S_0)\textrm{mes}(\Omega)}+\frac{1}{2}Ckq\sqrt{T_0}
< \frac{1}{2}\textrm{mes}(S_0).
$$
So, it follows from (5.5), (5.6) that:
$$
\|y^{(1)}-y^{(2)}\|_{(L^2(Q_{T_0}))^3}<\frac{1}{2}\|y^{(1)}-y^{(2)}\|_{(L^2(Q_{T_0}))^3}.
$$
Therefore $y^{(1)} \equiv y^{(2)}$ on $[0, T_0]$, and thus by (5.2),
$z_i^{(1)} \equiv z_{i}^{(2)}$ for $i = 1, \dots, n$ on $[0, T_0]$.
Contradiction. This ends the proof of Theorem 2.1.

\bigskip
{\bf 6. Conclusion.} In this paper we introduced   a new  hybrid model describing the locomotion of  a ``small'' swimmer  in the incompressible 3-$D$ fluid.  The model consists  of two  coupled systems of equations: one is a PDE  system--  for the fluid, governed by the 
nonstationary 3-$D$ Stokes  equation (typically associated with the low Reynolds numbers), and the other is an ODE  system-- for the position of  swimmer in it. It is  assumed that the swimmer's body  consists of 
finitely many subsequently connected parts, identified with the fluid they occupy,  which in turned are linked by the
 rotational  and elastic Hooke's forces. We investigated the well-posedness of this model. Namely, in suitable function spaces  we obtained  the  existence, uniqueness and regularity results for  its solutions. These results are to be used for the derivation of a formula for asymptotically small motions of this type of  3-$D$ swimmers, see \cite{Kh3}-\cite{Kh4} (and \cite{Kh2} for the 2-$D$ case).  They can also be instrumental  to study controllability properties of this swimming model \cite{Kh2}.  
 Models like this are of interest in biological
and engineering applications dealing with the study and design of self-propulsion systems in
fluids.

\end{document}